\documentclass[a4paper,
               12pt,
               pdftex,
               normalheadings,
               headsepline,
               footsepline,
               onecolumn,
               headinclude,
               footinclude,
               DIV14,
               abstracton]
{scrartcl}

\usepackage
{
    graphicx,
    amssymb,
    amsmath,
    amsthm,
    xcolor,
    dsfont,
    algpseudocode,
    authblk,
}

\usepackage{blindtext} 

\usepackage[bf]{caption}
\usepackage[colorlinks,
            pdffitwindow=false,
            plainpages=false,
            pdfpagelabels=true,
            pdfpagemode=UseOutlines,
            pdfpagelayout=SinglePage,
            bookmarks=false,
            colorlinks=true,
            hyperfootnotes=false,
            linkcolor=blue,
            citecolor=green!50!black]
{hyperref}

\usepackage{enumerate}
\usepackage{algorithm}

\usepackage{steinmetz}

\DeclareMathAlphabet{\mathpzc}{OT1}{pzc}{m}{it}

\newcommand{\R}{\mathbb{R}}
\newcommand{\C}{\mathbb{C}}
\newcommand{\Q}{\mathbb{Q}}
\newcommand{\N}{\mathbb{N}}
\newcommand{\Z}{\mathbb{Z}}
\newcommand{\mf}[1]{\mathpzc{#1}}

\providecommand{\abs}[1]{\left\lvert #1 \right\rvert}

\newcommand{\minsub}[1]{\hbox{\raise-2mm\hbox{$\textstyle \min \atop \scriptstyle {#1}$} }}

\newcommand\xqed[1]{\leavevmode\unskip\penalty9999 \hbox{}\nobreak\hfill \quad\hbox{#1}}

\newcommand{\littleo}{o \bigl( s^{-(1+\delta)} \bigr)}

\makeatletter
\renewcommand*\env@matrix[1][*\c@MaxMatrixCols c]{%
  \hskip -\arraycolsep
  \let\@ifnextchar\new@ifnextchar
  \array{#1}}
\makeatother

\begin{document}

\title{On BIBO stability of systems with irrational transfer function}
\author[*]{Ansgar Tr{\"a}chtler}
\affil[*]{\normalsize Heinz Nixdorf Institute, University of Paderborn, F{\"u}rstenallee~11, D-33098 Paderborn}

\maketitle

\begin{abstract}
We consider the input/output-stability of linear time-invariant single-input/single-output systems in terms of singularities of the transfer function F(s) in Laplace domain. The approach is based on complex analysis. A fairly general class of transfer functions F(s) is considered which can roughly be characterized by two properties: (1) F(s) is meromorphic with a finite number of poles in the open right half plane and (2), on the imaginary axes, F(s) may have at most a finite number of poles and branch points. For this class of systems, a complete and thorough characterization of BIBO stability in terms of the singularities in the closed right half-plane is developed where no necessity arises to differentiate between commensurate and incommensurate orders of branch points as often done in literature. A necessary and sufficient condition for BIBO stability is derived which mainly relies on the asymptotic expansion of the impulse response f(t). The second main result is the generalization of the well-known Nyquist citerion for testing closed-loop stability in terms of the open-loop frequency response locus to the class of systems under consideration which is gained by applying Cauchy's argument principle.
\end{abstract}

\pagestyle{myheadings}
\thispagestyle{plain}

\section{Introduction}
\label{sec:Introduction}
Non rational transfer functions arise in various context: \textit{spatially distributed systems} described by partial differential equations \cite{schober_1981_at,Fra_87,Deutscher2012} or, more general, \textit{infinite-dimensional systems} \cite{curtain_95}, or systems with one or more dead-time elements, e.g. \textit{retarded quasi-polynomial (RQ) meromorphic functions} \cite{zitek2008meromorphic}, or in \textit{fractional calculus} \cite{petras2008stability,matignon1996stability,matignon1998stability}. The approach to these systems usually relies on the powerful methods of functional analysis, and stability questions are formulated in terms of spectral properties of operators between Banach or Hilbert spaces.

The present contribution focuses on an approach based on complex analysis for characterizing input/output-stability of linear time-invariant single-input/single-output systems in terms of singularities of the transfer function $F(s)$ in Laplace domain. A fairly general class of transfer functions is considered which can roughly be characterized by two properties: (1) $F(s)$ is meromorphic with a finite number of poles in the open right half plane $\Re s > 0$ and (2), on the imaginary axes, $F(s)$ may have at most a finite number of poles and branch points. A  correct description is given in section \ref{sec:Openloop}. 

For this class of systems, a complete and thorough characterization of BIBO stability in terms of the singularities in the closed right half-plane $\Re s \ge 0$ is developed where no necessity arises to differentiate between commensurate and incommensurate orders of branch points as often done in literature, cf. e.g. \cite{matignon1998stability,petras2008stability}. A necessary and sufficient condition for BIBO stability is derived which mainly relies on the asymptotic expansion of the impulse response $f(t)=\mathcal L^{-1}\{F(s)\}$ \cite{Doe_74e}. This result extends the stability conditions reported in \cite{matignon1998stability} and differs from that in \cite{petras2008stability}. The second main result is the generalization of the well-known Nyquist citerion for testing closed-loop stability in terms of the open-loop frequency response locus to the class of systems under consideration which is gained by applying Cauchy's argument principle.

The paper is organized as follows: in section \ref{sec:Stability} some stability definitions are reported together with their relation to the impulse response. Section \ref{sec:Singularities} gives an introduction to isolated singularities, especially branch points and analytic continuation which forms the important background for the subsequent sections and the complete contribution. In section \ref{sec:Openloop}, the assumptions on the open-loop systems under consideration are formulated and the fundamental stability criterion is extended to this class of systems. In section \ref{sec:Closedloop} it is shown that, under the same assumptions to the open-loop system, the extended criterion also holds for the closed-loop systems. In section \ref{sec:NyquistCrit} the Nyquist criterion for testing closed-loop stability is extended to the class of irrational transfer functions considered here, and in section \ref{sec:Conclusion} the results are summarized.

\section{Stability}
\label{sec:Stability}

We consider linear, time-invariant (LTI) SISO systems whose input-output behavior is described by the convolution integral
\[ y(t) = \int g(t-\tau) u(\tau) d\tau \, .\]
$u(t)$, $y(t)$ denote input and output signal, respectively, with $u(t) = y(t) = 0$ for $t<0$, and the convolution kernel $g(t)$ represents the impulse response of the system. We assume that the system is causal which implies that also $g(t) = 0$ for $t<0$. Then the convolution integral results in 
\begin{equation} 
y(t) = g(t) \ast u(t) = \int \limits_{\tau=0}^t g(t-\tau) u(\tau) d\tau \, .
\label{eq:faltung}
\end{equation}
If $g(t)$, $u(t)$ grow at most exponentially as $t \rightarrow \infty$, i.e. there exist constants $M, \alpha > 0$ such that $\abs{g(t)}, \abs{u(t)} \le M e^{-\alpha t}$, then the Laplace integrals
\[
G(s) = \mathcal L \{ g(t) \} = \int \limits_{t=0}^\infty g(t) e^{-s t} dt \, , \quad 
U(s) = \mathcal L \{ u(t) \}= \int \limits_{t=0}^\infty u(t) e^{-s t} dt
\]
exist and converge absolutely for $\Re s \ge \alpha$. $G(s)$, $U(s)$ are bounded there and holomorphic for $\Re s > \alpha$ allowing the convolution integral \eqref{eq:faltung} to be mapped to
$Y(s) = G(s) \cdot U(s)$ with $Y(s)$ absolutely converging for $\Re s \ge \alpha$ as well \cite{Doe_74e}.

For LTI systems, stability is most commonly defined as BIBO stability. A system is called \textit{BIBO stable} if $\abs{u(t)} < m_u$ implies $\abs{y(t)} < M_u$ with constants $m_u$, $M_u$ depending on $u(t)$. Sometimes also a stability definition is used considering the asymptotic behavior of the step response: a system is called \textit{`step response stable'}, in short \textit{SR stable}, if the system's step response 
\[
h(t) = \int \limits_{\tau=0}^t g(\tau) d\tau
\]
has a finite limit for $t \rightarrow \infty$. SR stability is equivalent to (simple) convergence of the integral of $g(t)$:
\begin{equation}
\text{SR stable} \quad \Leftrightarrow \quad 
  \Big| \int \limits_{t=0}^\infty g(t) dt \, \Big| \,  < \infty \, ,
\label{eq:sr-stab}
\end{equation}
whereas the stronger BIBO stability is equivalent to absolute convergence
\begin{equation}
\text{BIBO stable} \quad \Leftrightarrow \quad \int \limits_{t=0}^\infty \abs{g(t)} dt < \infty \, ,
\label{eq:bibo-stab}
\end{equation}
for details see, e.g. \cite{Foe13}, or other textbooks. In \cite{curtain_95}, exponential stability is considered. We call a system  \textit{$\beta$-exponentially stable} where $\beta \le 0$ is any non-positive real number if there exist positive constants $M$ and $a$ with $a > - \beta$ such that
\[
\abs{g(t)} \le M e^{-a t}  \, , \quad \forall t \ge 0 \,.
\]
If $\beta = 0$, we simply say the system is \textit{exponentially stable}. Again, there is an equivalent necessary condition in terms of an integral of the impulse response, which is also sufficient if $g(t)$ is bounded: 
\begin{align}
\beta \text{-exponentially stable} \quad & \Leftrightarrow  \quad
\int \limits_{t=0}^\infty \abs{g(t)} e^{-\beta_1 t} dt < \infty \quad \text{for some } \beta_1 < \beta \, , 
\label{eq:exp-stab} \\
 & \quad \qquad \mbox{and} \quad \abs{g(t)} \le M < \infty \, .
\nonumber
\end{align}
The proof is easy and given in the appendix. Thus, exponential stability ($\beta = 0$) implies BIBO stability. 

As an example, we consider a system with (bounded) impulse response $g(t) = (t+1)^{-3/2}$ which is clearly BIBO stable but not exponentially stable. It can be shown that the corresponding transfer function $G(s)$ has an asymptotic expansion $G(s) \sim 2 - 2 \sqrt{\pi s}$ as $s \rightarrow 0$ and thus has a branch point in $s = 0$. In the next example, let 
\[
g(t) = \left\{
\begin{array}{ll}
\frac{1}{\sqrt{t}}	  & 0 < t \le 1 \\ & \\
\frac{1}{\sqrt{t}^3}	&  t >  1 
\end{array}
\right.
\]
$g(t)$ is unbounded at $t=0$, yet absolutely integrable over $(0,\infty)$ and thus BIBO stable, but not exponentially stable either. The third example is 
\[
g(t) = \frac{2}{\pi} \frac{\sin(\omega_0 t)}{t} \quad \mbox{with} \quad
G(s) = \frac{2}{\pi} \arctan \left( \frac{\omega_0}{s} \right) = \frac{1}{j\pi} \log \frac{s+ j\omega_0}{s-j\omega_0} \, .
\] 
Here, $g(t)$ is integrable over $[0,\infty)$, yielding $\int_0^\infty g(t) dt = 1$ and thus SR stable, but not BIBO stable. The transfer function has logarithmic branch points at $s=\pm j \omega$.

For finite-dimensional systems, where the transfer function $G(s)$ is a proper rational function of $s$, these 3 stability definitions are equivalent, and the \textit{fundamental stability criterion} holds: \textit{Necessary and sufficient for stability is that all poles of $G(s)$ have negative real parts.} The fundamental stability criterion has early been extended to some classes of non-rational transfer functions, e.g. in \cite{Foe67} to SR-stability of closed-loop systems whose open-loop transfer functions consist of a strictly proper rational transfer function $G_o(s)$ and a dead-time element $e^{-Ts}$ with $T \ge 0$. 

Below, we will extend the fundamental stability criterion in the sense of BIBO stability to  a wider class of LTI systems with non-rational transfer functions and analyse the role of the singularities of $G(s)$ with special regard to branch points.

\section{Isolated Singularities}
\label{sec:Singularities}

In the literature on complex analysis, isolated singular points are not treated uniformly with respect to branch points. A frequently used definition is that given a function $F(s)$ which is holomorphic in some open region $D \subset \C$ except in $s_0 \in D$, then $s_0$ is called an isolated singular point of $F(s)$. According to this definition, there exist exactly three kinds of isolated singular points: removable singular points, poles, and essential singularities \cite{Rem_1_2002}, whereas branch points do not belong to isolated singularities.

\begin{figure}[htb]
	\centering
		\includegraphics[width=0.4\textwidth]{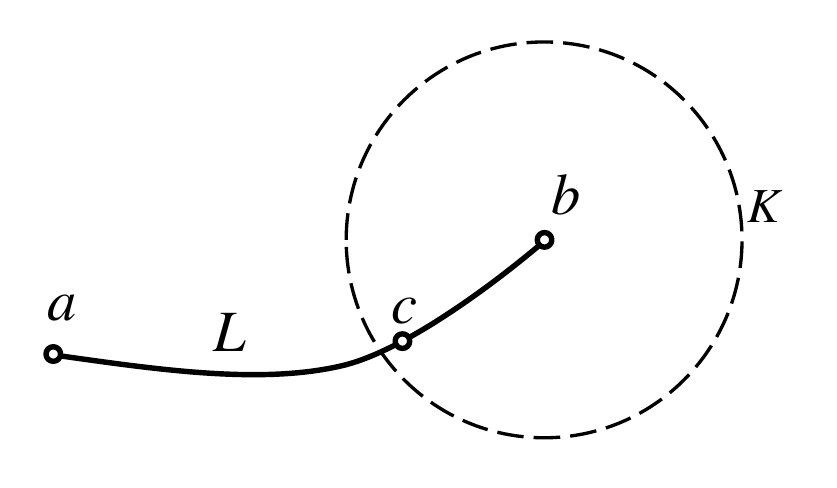}
	\caption{Isolated singularity ($b$) and analytic continuation from $a$ along $L$}
	\label{fig:fig_isol_sing}
\end{figure}

In our contribution, we follow W. I. Smirnov \cite{smirnov_eng_Vol_3-2_1964} who defines isolated singularities via the analytic continuation of holomorphic function elements. Starting from some point $a \in \C$ where $F(s)$ is holomorphic, this function element is extended along a curve $L$ (cf. fig. \ref{fig:fig_isol_sing}). If the analytic continuation is possible up to some point $b$ exclusively and not beyond, then $b$ is a singular point of $F(s)$ for the analytic continuation along $L$.

Now a circle $K$ is drawn around $b$ and the section $cb$ of $L$ lying inside $K$ is considered. If a circle $K$ exists such that the function elements starting in the points of section $cb$ can be analytically continued along every curve inside $K$, but not hitting $b$ then $b$ is called an \textit{isolated singular point} of $F(s)$ for the curve $L$.

Either the analytic continuation along all curves in $K$ with the same start and end point yield the same function element - then we call $F(s)$ single-valued or unambiguous - or there will result different function elements - then we call $F(s)$ multi-valued or ambiguous. If $F(s)$ is single-valued, then it is holomorphic in the whole interior of $K$ except $b$ and can be developed in a  series around $b$. The point $b$ is - with respect to the analytic continuation along $L$ - either a removable singularity, or a pole, or an essential singularity.

If the analytic continuations in $K$ are multi-valued, then $b$ is a branch point, and either a finite number $m$ or infinitely many different function elements exist. In the former case, $b$ is called branch point of order $m-1$, in the latter $b$ is a branch point of infinite order. Subsequently, we will term branch points of finite order \textit{algebraic branch points}.

For algebraic branch points of order $m-1$, Smirnov introduces a Riemann surface where $F(s)$ can be unambiguously defined (possibly except $b$) and  developed in a \textit{Laurent series} with respect to the argument $\sqrt[m]{s-b}$ yielding 
\begin{equation}
F(s) = \sum_{\nu=-\infty}^{\infty} c_\nu \left( \sqrt[m]{s-b} \right)^\nu 
     = \sum_{\nu=-\infty}^{\infty} c_\nu \left( s-b \right)^{\nu/m} \, .
\label{eq:Laurent_series}
\end{equation}
If the series \eqref{eq:Laurent_series} does not contain any summands with negative $\nu$,
\[
F(s) = c_0 + c_1 \sqrt[m]{s-b} + c_2 \sqrt[m]{s-b}^{\, 2} + \ldots \, ,
\]
this kind of algebraic branch point is called "regular-type branch point" and it holds $F(s=b)=c_0$. 
If the number of negative exponents in \eqref{eq:Laurent_series} is finite, the branch point is called "polar-type branch point". Finally, if the number of negative exponents is infinite, this is called "essentially-singular-type branch point", cf. \cite{smirnov_eng_Vol_3-2_1964}.

We extend our considerations to a more general class of functions including also infinite-order branch points of the type $(s-b)^\kappa$ with $\kappa \in \R$. Instead of a Laurent series, the function $F(s)$ is now expressed near $b$ by an \textit{asymptotic expansion}. According to the notation used in \cite{Doe_74e}, we develop $F(s)$ in a sequence of comparison functions $c_\nu (s-b)^{\kappa_\nu}$, $\nu=0,1,\ldots$ which represent $F(s)$ with increasing accuracy:
\begin{equation}
F(s) = \sum_{\nu=0}^{n} c_\nu \left( s-b \right)^{\kappa_\nu} + o\bigl( (s-b)^{\kappa_n}\bigr)  \quad \mbox{as} \quad s \rightarrow b \, .
\label{eq:asymptotic_expansion1}
\end{equation}
Herein, $\kappa_0 < \kappa_1 < \ldots$ and $c_\nu \neq 0$. If equation \eqref{eq:asymptotic_expansion1} is valid for each $n=0,1,\ldots$, $F(s)$ has an asymptotic expansion which is symbolically expressed by "$\approx$":
\begin{equation}
F(s) \approx \sum_{\nu=0}^{\infty} c_\nu \left( s-b \right)^{\kappa_\nu} \quad \mbox{as} \quad s \rightarrow b\, , \quad \mbox{with} \quad \kappa_0 < \kappa_1 < \ldots \; .
\label{eq:asymptotic_expansion2}
\end{equation}
Analogous to algebraic branch points, if $\kappa_0 \ge 0$, i.e. equation \eqref{eq:asymptotic_expansion2} does not contain any negative exponents, we call $b$ a \textit{regular branch point}, whereas for $\kappa_0 < 0$ we call $b$ a \textit{polar branch point}. Branch points with an infinite number of negative exponents are excluded from our consideration. 

We give some examples for functions with different kinds of branch points (all in $s=0$):

$ $ 

Algebraic branch points of order $m-1$ ($m\in \N$)
\begin{enumerate}
	 \item regular-type branch points: $\displaystyle \sqrt[m]{s}$, $\displaystyle e^{\sqrt[m]{s}}$, $\displaystyle \frac{1}{\sqrt{s}+2}$
	
The last example has a branch point of order $1$ in $s=0$ and a pole in $s=4$ lying in that leaf of the Riemann surface where $\sqrt 4 = -2$. In the other leaf where $\sqrt 4 = 2$, the point $s=4$ is regular.
	 \item polar-type branch points: $\displaystyle \frac{1}{\sqrt[m]{s}}$
	 \item essentially-singular-type branch point: $\displaystyle e^{1/{\sqrt[m]{s}}}$
\end{enumerate}
	
Branch points of infinite order
\begin{enumerate} 
\setcounter{enumi}{3}
  \item $\log(s)$
  \item $s^\alpha$ with $\alpha \in \R\backslash\Q$
\end{enumerate}
	
Removable branch point
\begin{enumerate}
\setcounter{enumi}{5}
  \item $\displaystyle \frac{\sinh(\sqrt{s}\cdot x)}{\sinh(\sqrt{s})}\, ,\quad 0\le x\le 1$
\end{enumerate}

Transfer functions of this kind arise in the context of the heat equation \cite{Fra_87, curtain_95}. That $s=0$ is indeed a removable branch point can be seen by expanding $\sinh$ in a power series:
\[ \frac{\sinh(\sqrt{s} \cdot x)}{\sinh(\sqrt{s})} 
= \frac{ \sum \limits_{k=0}^{\infty} \displaystyle \frac{ (\sqrt{s} \cdot x)^{2k+1}}{(2k+1)!} }
{ \sum \limits_{k=0}^{\infty} \displaystyle \frac{ (\sqrt{s} )^{2k+1}}{(2k+1)!} }
 = \frac{ \sqrt{s} \cdot x \sum \limits_{k=0}^{\infty} \displaystyle  \frac{ (\sqrt{s} \cdot x)^{2k}}{(2k+1)!} }
{ \sqrt{s}  \sum \limits_{k=0}^{\infty} \displaystyle \frac{ (\sqrt{s} )^{2k}}{(2k+1)!} }
 = x \cdot \frac{ \sum \limits_{k=0}^{\infty} \displaystyle \frac{ (s \cdot x^2 )^k}{(2k+1)!} }
{ \sum \limits_{k=0}^{\infty} \displaystyle \frac{ s^k}{(2k+1)!} } \, .
\]
We conclude this section with some remarks:

\begin{enumerate}
\item 
From a mathematical point of view, if $s=a$ is a removable singularity of $F(s)$, it is reasonable to define $F(a) = \lim_{s\rightarrow a} F(s)$; by this, $F(s)$ becomes regular in $s=a$, like in the last example. In context of control issues, however, this could stand for cancellation of singularities and should be handled with care. Nevertheless, below we will exclude removable singularities from our consideration. 

\item
Subsequently, we also exclude essential and logarithmic singularities and restrict ourselves to functions of the type \eqref{eq:asymptotic_expansion2}. The asymptotic expansion \eqref{eq:asymptotic_expansion2} can also be used to represent $F(s)$ in the neighborhood of algebraic branch points of regular and polar type ($\kappa_\nu \in \Q$) or poles ($\kappa_\nu \in \Z$), then yielding the respective Laurent series. 

\item
In the presence of branch points, functions are multi-valued. Unambiguity is achieved on a Riemann surface with several (finite or infinite) leaves. In general, the existence of singularities depends on the path taken for analytic continuation, cf. example $(\sqrt{s}+2)^{-1}$ above, and one might ask which leaf has to be chosen. In our context, however, it will always be uniquely determined which leaf has to be taken since $F(s)$ represents the transfer function of a technical system. Applying Laplace transform to the (uniquely defined) impulse response $f(t)$ yields unique $F(s)$ in some right half plane where analytic continuation can start from.
\end{enumerate}


\section{Open-loop stability}
\label{sec:Openloop}

\begin{figure}[htb]
	\centering
		\includegraphics[width=0.75\textwidth]{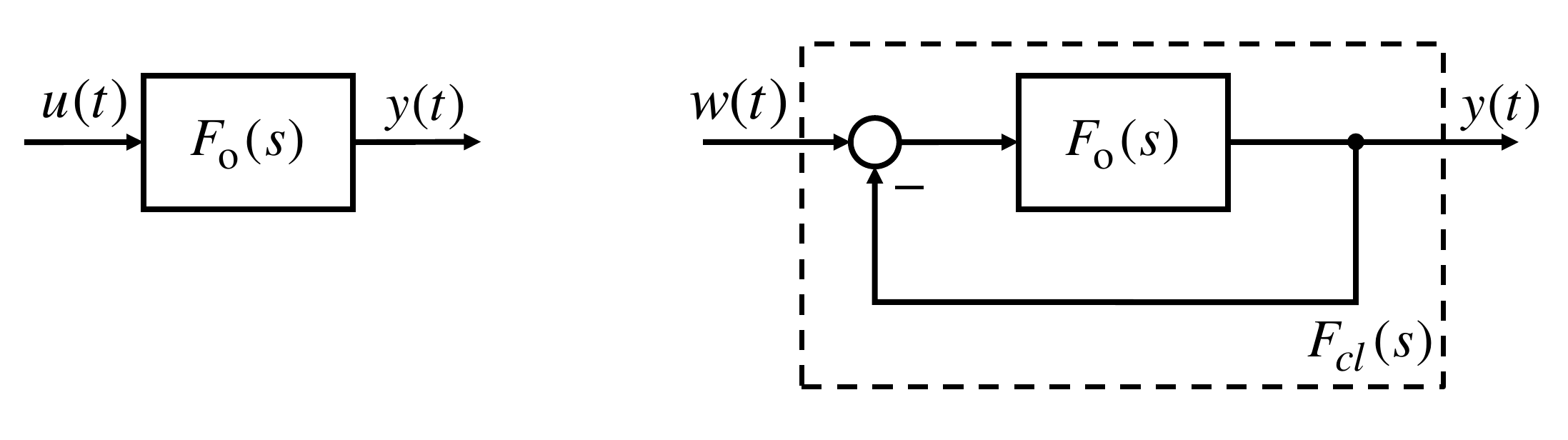}
	\caption{Open-loop and closed-loop system}
	\label{fig:Open_and_closed_loop}
\end{figure}

We consider the open-loop system shown in fig. \ref{fig:Open_and_closed_loop} on the left. $f_o(t)$ denotes the impulse response and $F_o(s) = \mathcal L \{ f_o(t) \}$ the corresponding transfer function. The aim is to formulate a necessary and sufficient condition for BIBO stability in terms of singularities of $F_o(s)$. Since  closed-loop stability will be considered in the next section, we aim to characterize the class of transfer functions under consideration in a way that also the closed-loop transfer function $F_{cl}(s) = \frac{F_o(s)}{1+F_o(s)}$ (cf. fig. \ref{fig:Open_and_closed_loop} on the right) belongs to this class. 

For the open-loop system, we meet the following assumptions (A1 ... A3).
\begin{description}
\item[\textnormal{A1: }] 
There is some $\alpha \ge 0$ such that
\begin{equation}
	\int \limits_{t=0}^\infty \abs{f_o(t)} e^{-\alpha t} dt < \infty \, .
	\label{eq:assump_1}
\end{equation}
\end{description}

Thus, $F_o(s) = \mathcal L \{ f_o(t) \}$ exists in the half-plane $\Re s \ge \alpha$ and is absolutely convergent there permitting the description of the transfer behavior in the frequency domain. Furthermore, $F_o(s)$ is bounded for $\Re s \ge \alpha$ and holomorphic for $\Re s > \alpha$ (\cite{Doe_74e}). 

\begin{figure}[htb]
	\centering
		\includegraphics[width=0.8\textwidth]{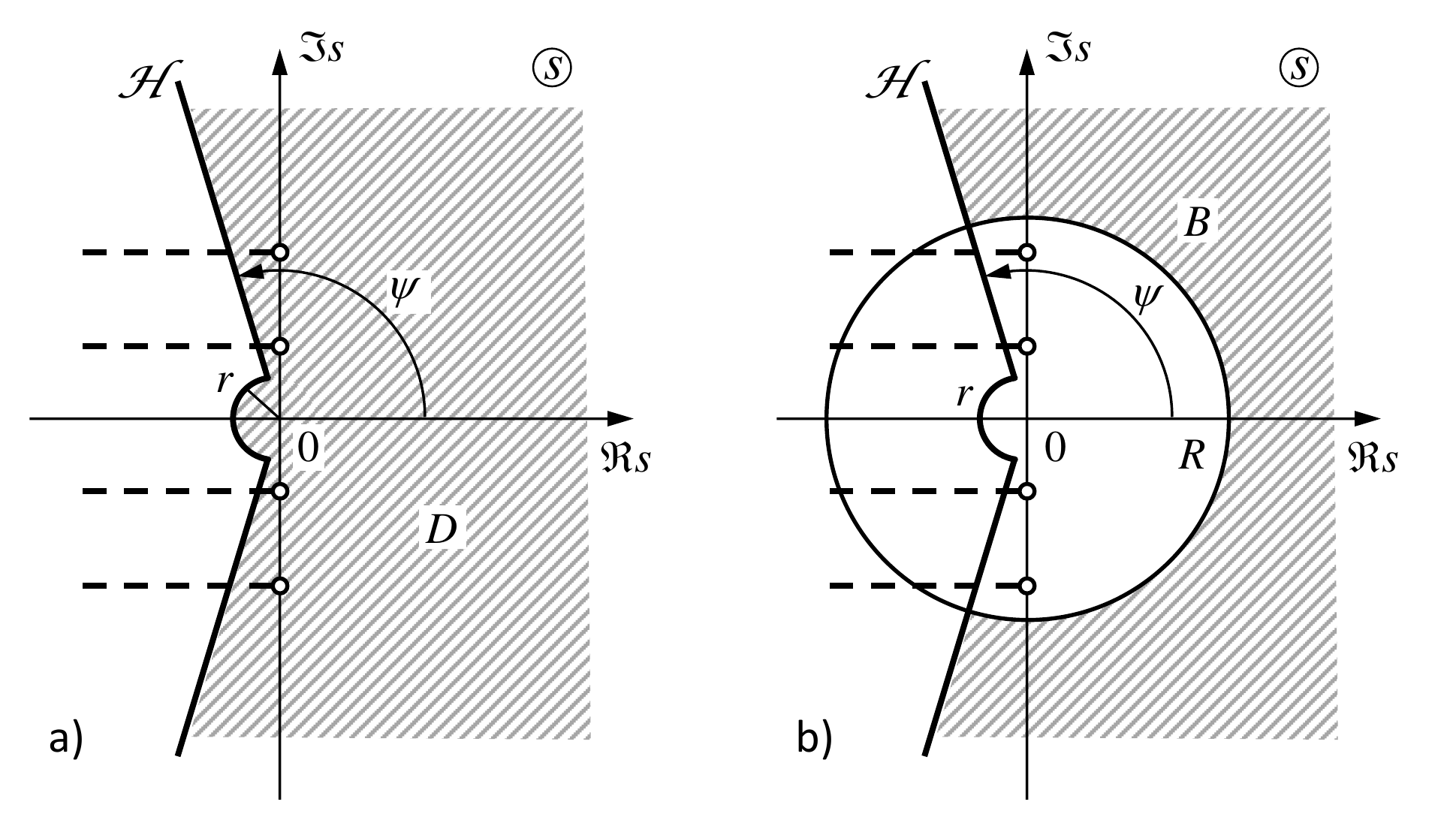}
	\caption{a) Region $D$ and location of possible branch points ("o") and cuts (- - -),
					\newline
	         b) Region $B$ of holomorphy of $F_o(s)$}
	\label{fig:fig_D_psi}
\end{figure}

Now, an open region $D$ is defined which is bounded to the left by a "hooked angle" $\mf{H} $ with center in $s=0$, circular arc with radius $r > 0$ and half opening angle $\psi$ ($\frac{\pi}{2}<\psi < \pi$) as shown in fig. \ref{fig:fig_D_psi} a).

\begin{description}
\item[\textnormal{A2: }] 
In $D$, $F_o(s)$ is assumed to be meromorphic with possible exception of a finite number of regular or polar branch points $b_i$ on the imaginary axis together with their corresponding cuts (parallels to the negative real axis). In the neighborhood of the branch points, $F_o(s)$ has the asymptotic expansion 
\begin{eqnarray}
F_o(s) &\approx &\sum_{\nu=0}^{\infty} c_{i\,\nu} \left( s-b_i \right)^{\kappa_{i\,\nu}} \quad \mbox{as} \quad s \rightarrow b_i\, , 
	\label{eq:assump_2}
\\
& & \qquad \quad \mbox{with} \quad c_{i\,\nu} \ne 0 \, , \quad \kappa_{i\,0} < \kappa_{i\,1} < \ldots \, ,\quad i=1,2,\ldots \; . \nonumber
\end{eqnarray}
Moreover, the number of negative exponents $\kappa_{i\,\nu} <0$ shall be finite.

\item[\textnormal{A3: }] 
In $D$, $F_o(s)$  shall asymptotically behave like 
\begin{equation}
F_o(s) = \frac{K}{s} + o \bigl( s^{-(1+\delta)} \bigr) =  \frac{K}{s} \bigl( 1 + o(s^{-\delta}) \bigr) \quad \mbox{for } s\rightarrow \infty \, , 	
\label{eq:assump_3}
\end{equation}
with some $\delta >0$. 
\end{description}

A2 implies that in every bounded region within $D$, $F_o(s)$ has at most a finite number of  branch points and cuts, as well as a finite number of poles, but no further singular points. In general, meromorphic functions can have poles which accumulate at infinity. But this is excluded by A3: Since $F_o(s) \rightarrow 0$ for $s\rightarrow \infty$ in $D$, there exists a radius $R>0$ such that $F_o(s)$ is holomorphic for $\abs{s} > R$,  $s\in D$. This region of holomorphy is denoted by $B$ and depicted in fig. \ref{fig:fig_D_psi} b) as the shaded area. Thus, by A2 and A3, $F_o(s)$ has at most a finite number of branch points and cuts, and a finite number of poles within entire $D$.

Now, we give a {\em generalization of the fundamental stability criterion} for systems satisfying the assumptions A1, A2, A3:
\begin{equation}
\begin{minipage}[c]{0.75\textwidth} \em
A linear time-invariant system of convolution type \eqref{eq:faltung} with impulse response $f_o(t)$ satisfying assumptions A1, A2, A3 is BIBO stable if and only if all singular points of $F_o(s)$ lie in the left half plane $\Re s < 0$ with possible exception of a finite number of regular branch points on the imaginary axis.
\end{minipage} \qquad \qquad
\label{Satz:StabKrit}
\end{equation}

{\bf Proof:} 1. Necessity. Assume $F_o(s)$ is BIBO stable. Then due to \eqref{eq:bibo-stab}, $f_o(t)$ is absolutely integrable and $F_o(s)$ is holomorphic in $\Re s >0$ and bounded in $\Re s \ge 0$ \cite{Doe_74e}. Thus the closed half plane $\Re s \ge 0$ is devoid of poles. On the imaginary axis, there could be branch points with $0 \le \kappa_{i\,0} < \kappa_{i\,1} < \ldots$, i.e. regular branch points.

2. Sufficiency. For proving the sufficient part, we assume that all singular points of $F_o(s)$ lie in the left half plane $\Re s < 0$ whereas on the imaginary axis we admit a finite number of regular branch points $b_i = j \omega_i$, $i=1,2,\ldots$ with asymptotic expansions according to \eqref{eq:assump_2} and $\kappa_{i\, \nu} \ge 0$. We have to show that this implies $f_o(t)$ to be absolutely integrable.

\begin{figure}[htb]
	\centering
		\includegraphics[width=0.4\textwidth]{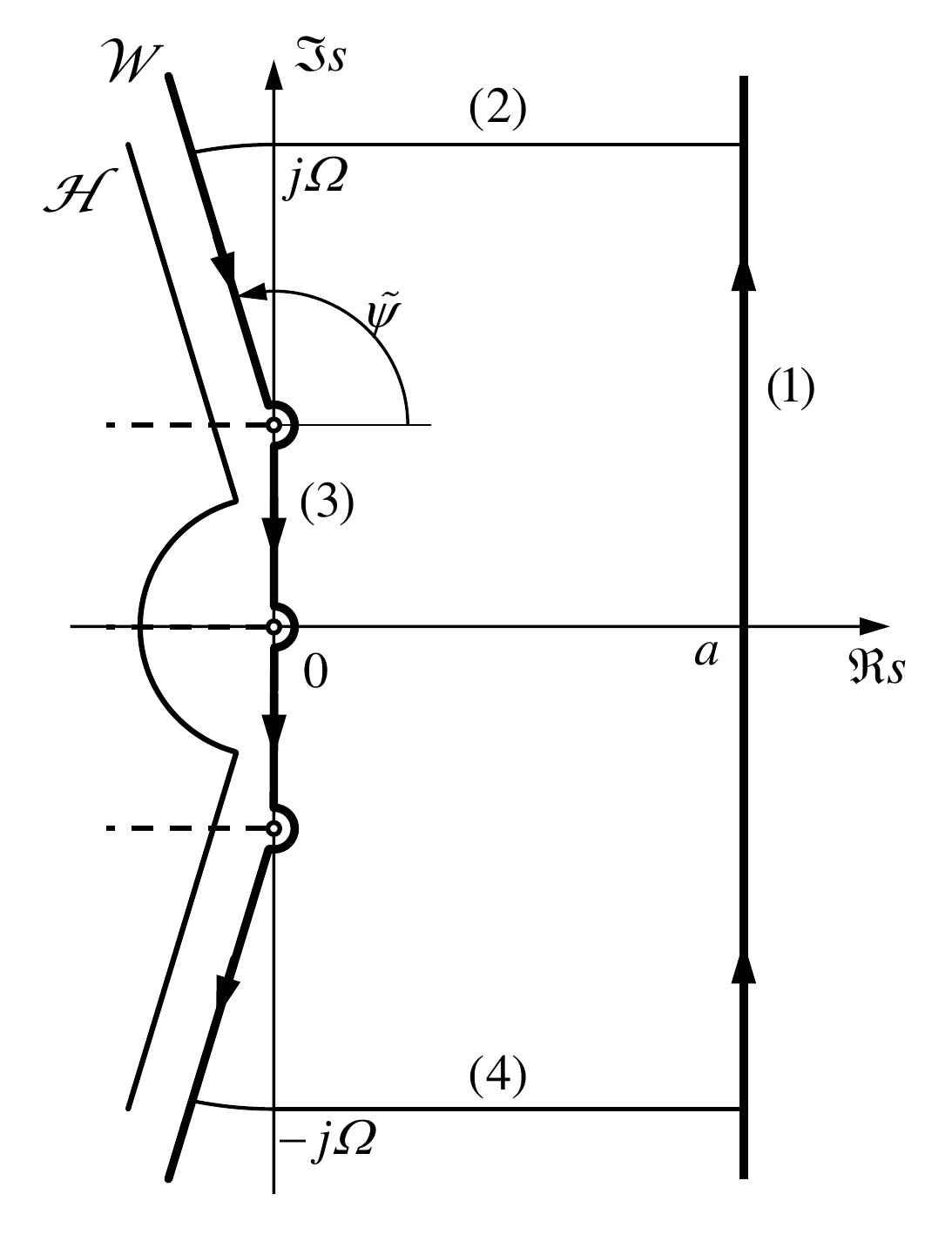}
	\caption{path of integration}
	\label{fig:fig_Umlaufintegral}
\end{figure}

We choose the closed contour shown in fig. \ref{fig:fig_Umlaufintegral} which is completely contained in region $D$ from fig. \ref{fig:fig_D_psi}. The contour consists of 4 parts: 

\begin{description}
\item[\textnormal{(1) }] 
$s=a+j\,y$, $ds=j\,dy$, $-\Omega \le y \le \Omega$, ($a \ge \alpha$ lying in the half plane of absolute convergence of $F_0(s)$, cf. A1),
\item[\textnormal{(2) }] 
$s=x+j\,\Omega$, $ds=dx$, $a \ge x\ge 0$, and $s=\Omega\,e^{j \,\varphi}$, $ds=j\,s\,d\varphi$, $\pi/2 \le \varphi \le \tilde{\psi}$, ($\Omega$ sufficiently large), 
\item[\textnormal{(3) }] 
the curve $\mf{W}$ of inclined rays ($\pi/2 < \tilde{\psi} \le \psi$, small arcs circumventing branch points to the right and parts of imaginary axis in between,
\item[\textnormal{(4) }] 
$s=\Omega \,e^{j \,\varphi}$, $ds=j\,s\,d\varphi$, $-\tilde{\psi} \le \varphi \le -\pi/2$, and 
$s=x-j\Omega$, $ds=dx$, $0 \le  x \le a$.
\end{description}

Since, by assumption, $F_o(s)$ is holomorphic on the contour and inside, the integral of $F_o(s) \, e^{s\,t}$ where $t>0$ yields
\[ \int \limits_{(1)+(2)+(3)+(4)} F_o(s) \, e^{s\,t} ds = 0 \, .
\]
For $\Omega \rightarrow \infty$, the integrals along (2) and (4) tend to zero due to assumption A3, cf. \cite{Doe_74e}, thus the integrals along (1) and (3) sum up to zero. Using the inversion formula of Laplace transform, we obtain
\begin{equation}
f_o(t) = \frac{1}{2 \pi j} \int \limits_{s=a-j\infty}^{a^+j\infty} F_o(s) \, e^{s\,t} ds 
=  - \int \limits_{\mf{W}} F_o(s) \, e^{s\,t} ds \, , \quad t>0 \, .
\end{equation}

For the integral on the right, theorem 37.2 in \cite{Doe_74e} applies. According to this theorem, the asymptotic behaviour of the original function $f_o(t)$ as $t \rightarrow \infty$ is represented by superposition of the asymptotic expansions
\begin{equation}
\sum_{\nu=0}^{\infty} \frac{c_{i\,\nu}}{\Gamma(-\kappa_{i\,\nu})} \, 
\frac{ e^{j\omega_i t}}{t^{\kappa_{i\,\nu}+1}} \, , \quad i=1,2,\ldots \; .
\label{eq:asymptotic_expansion_time}
\end{equation}
Herein, $\Gamma(.)$ is the gamma function with the understanding that, for $\kappa_{i\,\nu} = 0, 1, 2\ldots$, 
\begin{equation}
\frac{1}{\Gamma(-\kappa_{i\,\nu})} = 0 \, .
\label{eq:gamma_function}
\end{equation}
The asymptotic behavior of $f_o(t)$ is dominated by the term with smallest exponent $\kappa_{i\,\nu}$ contributing to \eqref{eq:asymptotic_expansion_time}. Due to \eqref{eq:gamma_function}, nonnegative integer values of $\kappa_{i\,\nu}$ have to be excluded since they do not contribute to \eqref{eq:asymptotic_expansion_time}. Thus the dominant term is given by the minimum value of all non-integer exponents: $\kappa^* = \minsub{i,\,\nu} \left\{ \kappa_{i\,\nu}\, | \,\kappa_{i\,\nu} \notin \N_0 \right\}$. Since by assumption, all $\kappa_{i\,\nu}$ are $\ge 0$, it holds that $\kappa^* >0$. Finally we obtain for the impulse response 
\begin{equation}
f_o(t) = \frac{1}{t^{\kappa^*+1}} \cdot O(1) \, , \quad \mbox{as } t \rightarrow \infty\,,
\label{eq:asymptotic_expansion_time2}
\end{equation}
and conclude that $f_o(t)$ is absolutely integrable over $(0,\infty)$ and thus BIBO stable. Hence the proof is complete. 

\xqed{$\blacksquare$}

\section{Closed-loop stability}
\label{sec:Closedloop}
In this section we consider BIBO stability of the closed-loop system depicted in fig. \ref{fig:Open_and_closed_loop} on the right and show that the stability criterion derived in the last section is also applicable to the closed-loop system. Therefore, we first show that the closed-loop behavior is also described by a transfer function in Laplace domain and secondly that the closed-loop transfer function satisfies the assumptions (A1) to (A3) if the open-loop transfer function does.

We start with defining the class $L_{1,\alpha}$ of functions which, when multiplied by $e^{-\alpha t}$, are Lebesgue integrable in the interval $(0,\infty)$ (with $\alpha \ge 0$). $L_{1,\alpha}$ is a Banach space under the norm
\[
\| f \|_{1,\alpha} := \int \limits_0^{\infty} |f(t)| \, e^{-\alpha t} \, dt \, .
\]
For the open-loop system in fig. \ref{fig:Open_and_closed_loop} where 
$	y(t) = f_o(t) \ast u(t)$, condition \eqref{eq:assump_1} now reads $f_o \in L_{1,\alpha}$.
It is easily proven that for input signals $u \in L_{1,\alpha}$, also the output $y$ belongs to $L_{1,\alpha}$. Thus we can interprete the open-loop system as a linear operator $A: L_{1,\alpha} \rightarrow L_{1,\alpha}$ which is bounded under the norm $\| A\| = \|f_o\|_{1,\alpha}$. 

Now we turn to the closed-loop system shown in fig. \ref{fig:Open_and_closed_loop} on the right. The time-domain representation of the input/output behavior reads
\[
y(t) = f_o(t) \ast \bigl( w(t) - y(t) \bigr)\,\quad \mbox{or}\quad y = A\, (w-y) \, .
\]
With the identity operator $I$, this can be converted into
\[
\bigl( I + A \bigr) \, y = A\, w \, .
\]
Due to \eqref{eq:assump_3}, for $f_o(t)$ given, we can choose $\alpha$ sufficiently large such that $\|A\|=\|f_o\|_{1,\alpha}<1$. Then, the corresponding Neumann series $\sum_{\nu=0}^\infty (-A)^\nu$ is uniformly convergent and $I+A$ has a bounded inverse, cf. e.g. \cite{Heuser_Funktionalanalysis}. Thus,
\[
y = \bigl( I + A \bigr)^{-1} A\, w \, \in L_{1,\alpha} \quad \mbox{with} \quad 
\|( I + A )^{-1} A\| \le \frac{\|A\|}{1-\|A\|} \, .
\]
We can conclude that the Laplace integral $Y(s)=\int_{t=0}^\infty y(t) e^{-s t} dt$ is absolutely convergent for $\Re s \ge \alpha$ and the closed-loop system can be represented in Laplace domain by the transfer function
\[
F_{cl}(s) = \frac{F_o(s)}{1+F_o(s)} \, 
\]
and in time domain by the impulse response $f_{cl}(t)=\mathcal{L}^{-1} \{F_{cl}(s)\}$. Since $\|f_{cl}\|_{1,\alpha}= \|( I + A )^{-1} A\|$, the norm of $f_{cl}(t)$ is bounded, thus $\int_0^{\infty} |f_{cl}(t)| \, e^{-\alpha t} \, dt < \infty$ and by this we have shown that A1 also holds for the closed-loop system.

Next we have to prove, that the assumptions A2 and A3 apply to $F_{cl}(s)$. For proving A3 we calculate
\[ F_{cl}(s)-\frac{K}{s} = \frac{F_o(s)}{1+F_o(s)}-\frac{K}{s} 
= \frac{F_o(s)-\frac{K}{s}\,(1+F_o(s))}{1+F_o(s)} \, .
\]
Replacing $F_o(s)$ by $\frac{K}{s} + \littleo$ in the numerator and by $o(1)$ in the denominator, respectively, yields
\[ 
F_{cl}(s)-\frac{K}{s} = \frac{\frac{K}{s}+\littleo-\frac{K}{s} -\frac{K}{s}\,(\frac{K}{s}+\littleo)}{1+o(1)} \, = \littleo 
\]
and finally we obtain
\[ 
F_{cl}(s) = \frac{K}{s} + \littleo \, .
\]

Next we turn to A2. $F_o(s)$ was assumed to be meromorphic in the open region $D$ with exception of branch points on the imaginary axis and the corresponding cuts. $F_{cl}(s)$ being the quotient of two meromorphic functions is itself meromorphic in the same region.

It remains to consider the behavior of $F_{cl}(s)$ in the neighborhood of branch points of $F_o(s)$. For simplicity, subsequently we will write $b$, $c_\nu$, $\kappa_\nu$ instead of $b_i$, $c_{i\,\nu}$, $\kappa_{i\,\nu}$. According to \eqref{eq:assump_2}, $F_o(s)$ is represented by its asymptotic expansion 
\begin{equation}
F_o(s) \approx \sum_{\nu=0}^{\infty} c_\nu (s-b)^{\kappa_\nu}	\quad \mbox{as} \quad s \rightarrow b
\label{eq:asymptotic_expansion4}
\end{equation}
where $\kappa_0 < \kappa_1 < \ldots$ and $c_\nu \ne 0$. If only the first term is interesting, we write short
\[
F_o(s) = c_0 (s-b)^{\kappa_0}	\bigl( 1+ o(1) \bigr) \, .
\]
Near $s=b$, the closed-loop transfer function behaves like
\begin{equation}
F_{cl}(s) = \frac{c_0 (s-b)^{\kappa_0}	\bigl( 1+ o(1) \bigr)}{1+c_0 (s-b)^{\kappa_0}	\bigl( 1+ o(1) \bigr)} \,. 
\label{eq:asymptotic_expansion_cl}
\end{equation}
It should be mentioned, that the "1" in the denominator of \eqref{eq:asymptotic_expansion_cl} has to be treated with care since it depends on $\kappa_0$ which term is dominant. Below, this will be illuminated in detail differentiating 4 cases with different values of $\kappa_0$ and $c_0$.

\begin{description}
\item[\textnormal{Case 1: }] $\kappa_0 < 0 \quad$ ($F_o(s)$ has a polar branch point in $s=b$)

$\qquad \quad \displaystyle
F_{cl}(s) = \frac{c_0 (s-b)^{\kappa_0}	\bigl( 1+ o(1) \bigr)}{c_0 (s-b)^{\kappa_0}	\bigl( 1+ o(1) \bigr)} = 1+ o(1) \, .
$

\item[\textnormal{Case 2: }] $\kappa_0 = 0$ and $c_0 \ne -1 \quad$ ($F_o(s)$ has a regular branch point in $s=b$)

$\qquad \quad \displaystyle
F_{cl}(s) = \frac{c_0 \bigl( 1+ o(1) \bigr)}{(1+c_0) 	\bigl( 1+ o(1) \bigr)} 
= \frac{c_0}{1+c_0} \, \bigl( 1+ o(1) \bigr) \, .
$

\item[\textnormal{Case 3: }] $\kappa_0 = 0$ and $c_0 = -1 \quad $ ($F_o(s)$ has a regular branch point in $s=b$)

Here we need the first two terms in $F_o(s)$:

$\qquad \quad \displaystyle
F_o(s) = c_0 (s-b)^{\kappa_0} + c_1 (s-b)^{\kappa_1}	\bigl( 1+ o(1) \bigr) 
 = -1 + c_1 (s-b)^{\kappa_1}	\bigl( 1+ o(1) \bigr) \,
$

where $\kappa_1 > \kappa_0 = 0$. For $F_{cl}(s)$ we obtain

$\qquad \quad \displaystyle
F_{cl}(s) = \frac{-1 + c_1 (s-b)^{\kappa_1}	\bigl( 1+ o(1) \bigr)}{1-1 + c_1 (s-b)^{\kappa_1}	\bigl( 1+ o(1) \bigr)} 
= -\frac{1}{c_1} \, (s-b)^{-\kappa_1} \bigl( 1+ o(1) \bigr) \, .
$

\item[\textnormal{Case 4: }] $\kappa_0 > 0 \quad$ ($F_o(s)$ has a regular branch point in $s=b$)

$\qquad \quad \displaystyle
F_{cl}(s) = \frac{c_0 (s-b)^{\kappa_0} \bigl( 1+ o(1) \bigr)}{1+c_0 (s-b)^{\kappa_0} \bigl( 1+ o(1) \bigr)} = c_0 (s-b)^{\kappa_0} \bigl( 1+ o(1) \bigr) \, .
$

\end{description}

In case 1, $F_o(s)$ has a polar branch point, and in cases 2, 3, and 4 regular branch points whereas $F_{cl}(s)$ has a polar branch point in case 3 and regular ones else. Summarized we have shown that the closed-loop transfer function belongs to the same class of transfer functions satisfying the assumptions A1, A2, A3 and by this we can conclude that also the stability criterion \eqref{Satz:StabKrit} applies to the closed-loop system.

\section{Nyquist criterion}
\label{sec:NyquistCrit}
\begin{figure}[htb]
	\centering
		\includegraphics[width=0.85\textwidth]{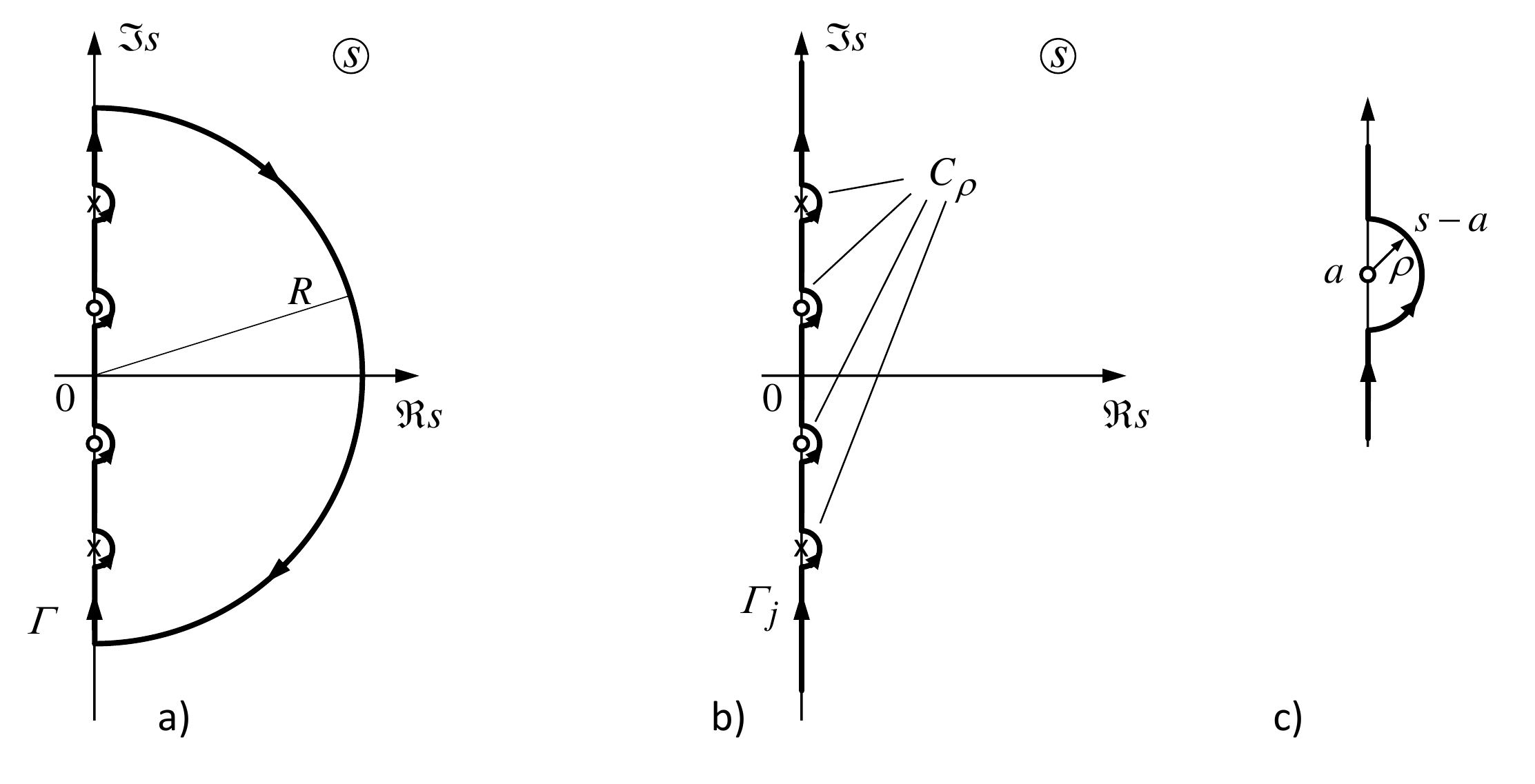}
	\caption{Nyquist contour}
	\label{fig:fig_Nyquist-contour}
\end{figure}

In this section, we extend the well-known Nyquist criterion to the class of open-loop transfer functions characterized by A1, A2, A3 from section \ref{sec:Openloop}. As usual, the fundament for the Nyquist criterion is Cauchy's argument principle which can be found in most textbooks on complex analysis, e.g. \cite{smirnov_eng_Vol_3-2_1964}. First, we define the closed contour $\Gamma$ shown in fig. \ref{fig:fig_Nyquist-contour} a). 

It consists of a semi-circle with (large) radius $R$ and the section of the imaginary axis between $-jR$ and $jR$. If the open-loop or closed-loop transfer functions have singular points (poles or branch points) on the imaginary axis, these are circumvented by small semi-circles $C_\varrho$ with radius $\varrho$. The radius $R$ is chosen sufficiently large such that the outer region $\abs{s} \ge R\, , \; \Re s \ge 0$ does not contain poles or branchpoints of $F_o(s)$ or $F_{cl}(s)$. $\varrho$ is chosen sufficiently small such that the interiors of the semi-circles, except their centers, do not contain further singularities. In the interior of $\Gamma$, $F_o(s)$ and $F_{cl}(s)$ are meromorphic. Then we can apply Cauchy's argument principle:
\begin{equation}
\int \limits_\Gamma \frac{F_o'(s)}{1+F_o(s)}ds = 2\pi \,j \,\bigl( P - Z \bigr)
\label{eq:Cauchy_Principle}
\end{equation}
 where $P$ and $Z$ denote the numbers of poles and zeros of $1+F_o(s)$ in the right half-plane $\Re s > 0$ (RHP), respectively (multiplicities are to be accounted for). (Since we integrate clockwise, the sign on the right-hand side of \eqref{eq:Cauchy_Principle} differs from the usual formulation).

The poles of $1+F_o(s)$ coincide with the poles of $F_o(s)$, whereas the zeros of $1+F_o(s)$ are the poles of the closed-loop transfer function $F_{cl}(s)$. We denote the number of RHP poles of $F_o$ by $P_{o,+}$ and the number of RHP poles of $F_{cl}$ by $P_{cl,+}$. Then, equation \eqref{eq:Cauchy_Principle} reads
\begin{equation}
\int \limits_\Gamma \frac{F_o'(s)}{1+F_o(s)}ds = 2\pi \,j \, \bigl( P_{o,+} - P_{cl,+} \bigr) \, .
\label{eq:Cauchy_Principle2}
\end{equation}

The integrand has the antiderivative 
$\log \bigl(1+F_o(s) \bigr) = \ln \abs{1+F_o(s)} + j  \phase{1+F_o (s)}$ and we obtain the value of the integral by determining the increase of the antiderivative when travelling along $\Gamma$. After one encirclement, the real part $\ln \abs{1+F_o(s)}$ ends in the same value when starting thus it doesn't contribute. Now let $R\rightarrow\infty$. Since $F_o(s) \rightarrow 0$ as $\abs{s}\rightarrow\infty$, on the semi-circle the argument $\phase{ 1+F_o (s)} \rightarrow 0$ and the integral has to be taken along the curve $\Gamma_j$ shown in fig. \ref{fig:fig_Nyquist-contour} b). If we denote the increase of argument along $\Gamma_j$ by the symbol $\Delta_{\Gamma_j}$, \eqref{eq:Cauchy_Principle2} becomes
\begin{equation}
\int \limits_{\Gamma_j} \frac{F_o'(s)}{1+F_o(s)}ds = j \, \underset{\Gamma_j}{\Delta} \, \phase{ 1+F_o (s) }= 2\pi \,j \, \bigl( P_{o,+} - P_{cl,+} \bigr) \, .
\label{eq:Cauchy_Principle3}
\end{equation}

In the last step, we consider the increase of argument along $\Gamma_j$ if the radii of the semi-circles tends to zero: $\varrho\rightarrow 0$. As shown in the appendix, for $s=a$ being a pole of $F_o(s)$ with multiplicity $k>0$, the increase of argument $\Delta_a$ is $ -k\,\pi$. If $s=a$ is pole of $F_{cl}(s)$ with multiplicity $k>0$, $\Delta_a = +k\,\pi$. If $s=a$ is a polar branch point of $F_o$ with dominant exponent $\kappa_0<0$, then $\Delta_a = \kappa_0\,\pi \, (<0)$.\footnote{ Thus, in a sense, $-\kappa_0$ can be interpreted as "multiplicity" of the polar branch point.}  If, finally, $s=a$ is a regular branch point of $F_o$ with $\kappa_0 = 0$, $c_0=-1$, the closed-loop systems has a polar branch point with dominant exponent $-\kappa_1 <0$ and the increase of argument is $\Delta_a = \kappa_1\,\pi$. In the remaining cases of imaginary branch points the increase of argument tends to zero as $\varrho \rightarrow 0$.

We combine these particular results by introducing the numbers $P_{o,j}$ and $P_{cl,j}$ of poles of $F_o(s)$ and of $F_{cl}(s)$ on the imaginary axis, respectively. For polar branch points on the imaginary axis, we define equivalent (nonnegative) numbers by \footnote{ In a similar manner $B_{o,j}$, $B_{cl,j}$ can be interpreted as "total number" of polar branch points including their "muliplicities".}
\begin{align}
B_{o,j} &:=  - \sum_{\stackrel{\scriptsize {\scriptstyle i} \mbox{ where}} 
{\kappa_{i,0}<0}} \kappa_{i,0} \qquad \mbox{for the polar branch points of $F_o$, and } 
\label{eq:branchnumber_open-loop} \\
B_{cl,j} &:=   \sum_{\stackrel{\scriptsize {\scriptstyle i} \mbox{ where}}
{\kappa_{i,0}=0\, , \; c_{i,0}=-1}} \hspace{-3mm} \kappa_{i,1} \quad \mbox{for the polar branch points of $F_{cl}$.} 
\label{eq:branchnumber_closed-loop}
\end{align}
Thus, equation \eqref{eq:Cauchy_Principle3} results in
\begin{align}
\underset{\Gamma_j}{\Delta} \, \phase{ 1+F_o (s) } &= 
\stackrel{+\infty}{\underset{\omega=-\infty}{\Delta^\bullet}} \phase{ 1 + F_o (j\omega)} 
- \pi \, (P_{o,j}+B_{o,j}) + \pi \, (P_{cl,j}+B_{cl,j}) \nonumber \\
&= 2\pi \, \bigl( P_{o,+} - P_{cl,+} \bigr) \, .\nonumber 
\end{align}
or
\begin{equation}
\stackrel{+\infty}{\underset{\omega=-\infty}{\Delta^\bullet}} \phase{ 1 + F_o (j\omega)} 
= \pi \, \bigl( 2 P_{o,+}  + P_{o,j}  + B_{o,j} \bigr) \, - \, 
  \pi \, \bigl( 2 P_{cl,+} + P_{cl,j} + B_{cl,j}\bigr) \, .
\label{eq:Winkeländerung_allgemein}
\end{equation}

\begin{figure}[htb]
	\centering
		\includegraphics[width=0.5\textwidth]{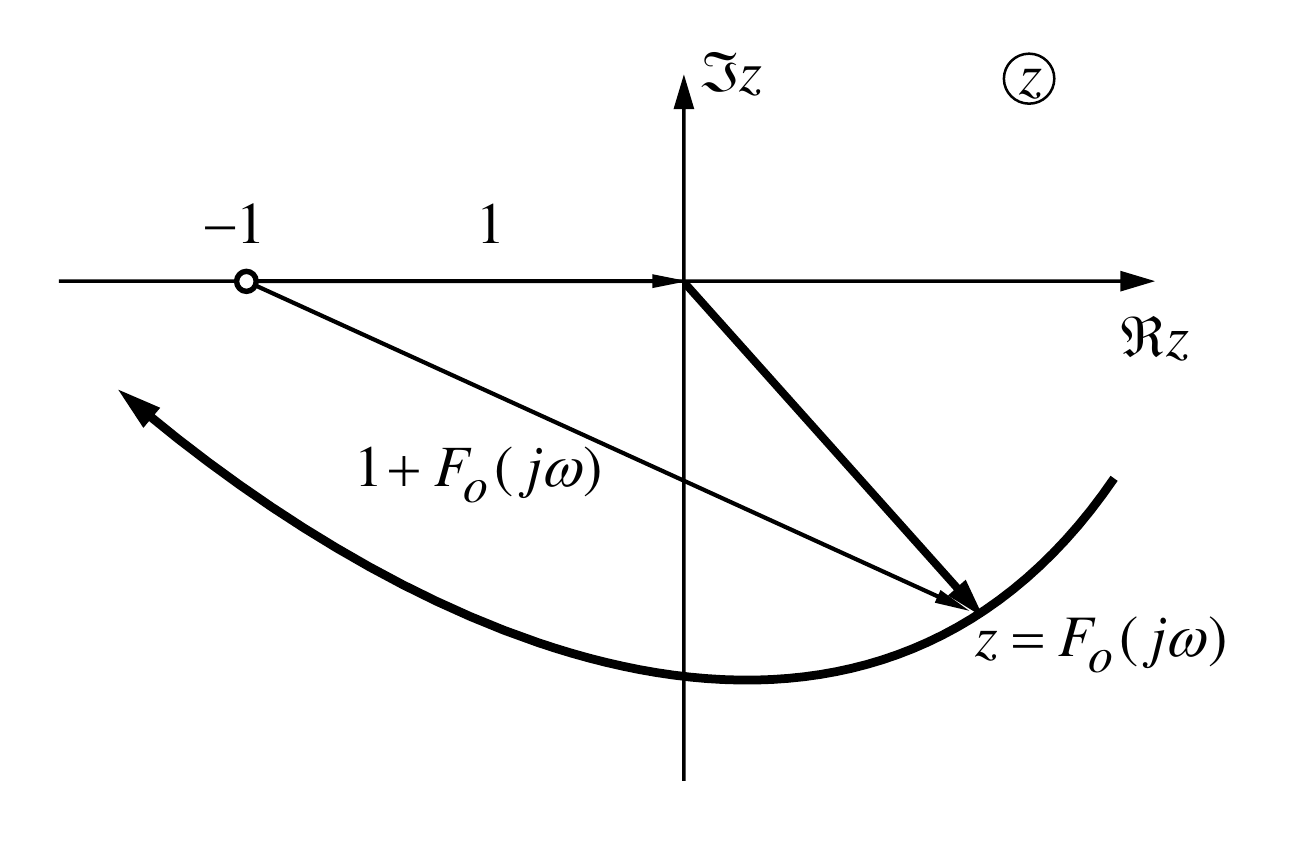}
	\caption{Ray from $z=-1$ to the current point on the Nyquist plot}
	\label{fig:fig_Fahrstrahl}
\end{figure}

Equation \eqref{eq:Winkeländerung_allgemein} admits an interpretation in terms of the Nyquist plot: $\stackrel{+\infty}{\underset{\omega=-\infty}{\Delta^\bullet}} \phase{ 1 + F_o (j\omega)}$ is the change in the angle of the ray pointing from the point $z=-1$ to the current point on the Nyquist plot $z=F_o(j\omega)$ (cf. fig. \ref{fig:fig_Fahrstrahl}) as $\omega$ varies on the \textit{punctured} $j$-axis (the singularities being removed) from $-\infty$ to $\infty$.

The relation \eqref{eq:Winkeländerung_allgemein} applies independently of any stability considerations. But due to the generalized stability criterion \eqref{Satz:StabKrit}, the closed-loop system is BIBO stable if and only if the numbers $P_{cl,+}$, $P_{cl,j}$, $B_{cl,j}$ are zero, i.e. the closed-loop systems has no poles in $\Re s \ge 0$ and no polar branch points on the imaginary axis which implies that 
$\stackrel{+\infty}{\underset{\omega=-\infty}{\Delta^\bullet}} \phase{ 1 + F_o (j\omega)} 
= \pi \, \bigl( 2 P_{o,+}  + P_{o,j}  + B_{o,j} \bigr)$. If, conversely, $\stackrel{+\infty}{\underset{\omega=-\infty}{\Delta^\bullet}} \phase{ 1 + F_o (j\omega)} 
= \pi \, \bigl( 2 P_{o,+}  + P_{o,j}  + B_{o,j} \bigr)$ holds, then equation \eqref{eq:Winkeländerung_allgemein} is only satisfied if $2 P_{cl,+} + P_{cl,j} + B_{cl,j} = 0$. Since $P_{cl,+}$, $P_{cl,j}$ and $B_{cl,j}$ are nonnegative, this implies $P_{cl,+} = P_{cl,j} = B_{cl,j} = 0$, thus the closed-loop system is BIBO stable.

Now we can formulate the \textit{generalized Nyquist citerion}:
\begin{equation}
\begin{minipage}[c]{0.75\textwidth} \em
Assuming the open-loop system to have a transfer function $F_o(s)$ satisfying the assumptions A1, A2, A3. The total number of poles of $F_o(s)$ in the right half plane and on the $j$-axis is $P_{o,+}$ and $P_{o,j}$, respectively. Furthermore, the "total number" of polar branch points on the $j$-axis according to \eqref{eq:branchnumber_open-loop} is $B_{o,j}$. Then, the closed-loop system is stable if and only if the change of argument of the ray pointing from $z=-1$ to the current point $z=F_o(j\,\omega)$ on the Nyquist plot, as $\omega$ varies on the punctured $j$-axis from $-\infty$ to $\infty$, equals
\[
\stackrel{+\infty}{\underset{\omega=-\infty}{\Delta^\bullet}} \phase{ 1 + F_o (j\omega)} 
= \pi \, \bigl( 2 P_{o,+}  + P_{o,j}  + B_{o,j} \bigr) \, .
\label{eq:Winkeländerung_stabil}
\]
\end{minipage} \qquad \qquad
\label{Satz:Nyquist}
\end{equation}

It should be noted that in literature when the open loop system has poles on the imaginary axis, the limit $\varrho \rightarrow 0$ is often not performed. Instead, semi-circles with small, but finite radius $\varrho$ to the \textit{left-hand side} of the singularities are considered (in that case, the Nyquist plot is a finite closed curve). For branch points on the imaginary axis, however, this procedure does not work, since integration along semi-circles on the left-hand side would cross the cuts and thus is not allowed. If, on the other hand, the change of argument along finite semi-circles on the right-hand side is considered, \eqref{eq:Cauchy_Principle3} reveals no information on potential singularities inside the semi-circle.

\begin{figure}[htb]
	\centering
		\includegraphics[width=0.5\textwidth]{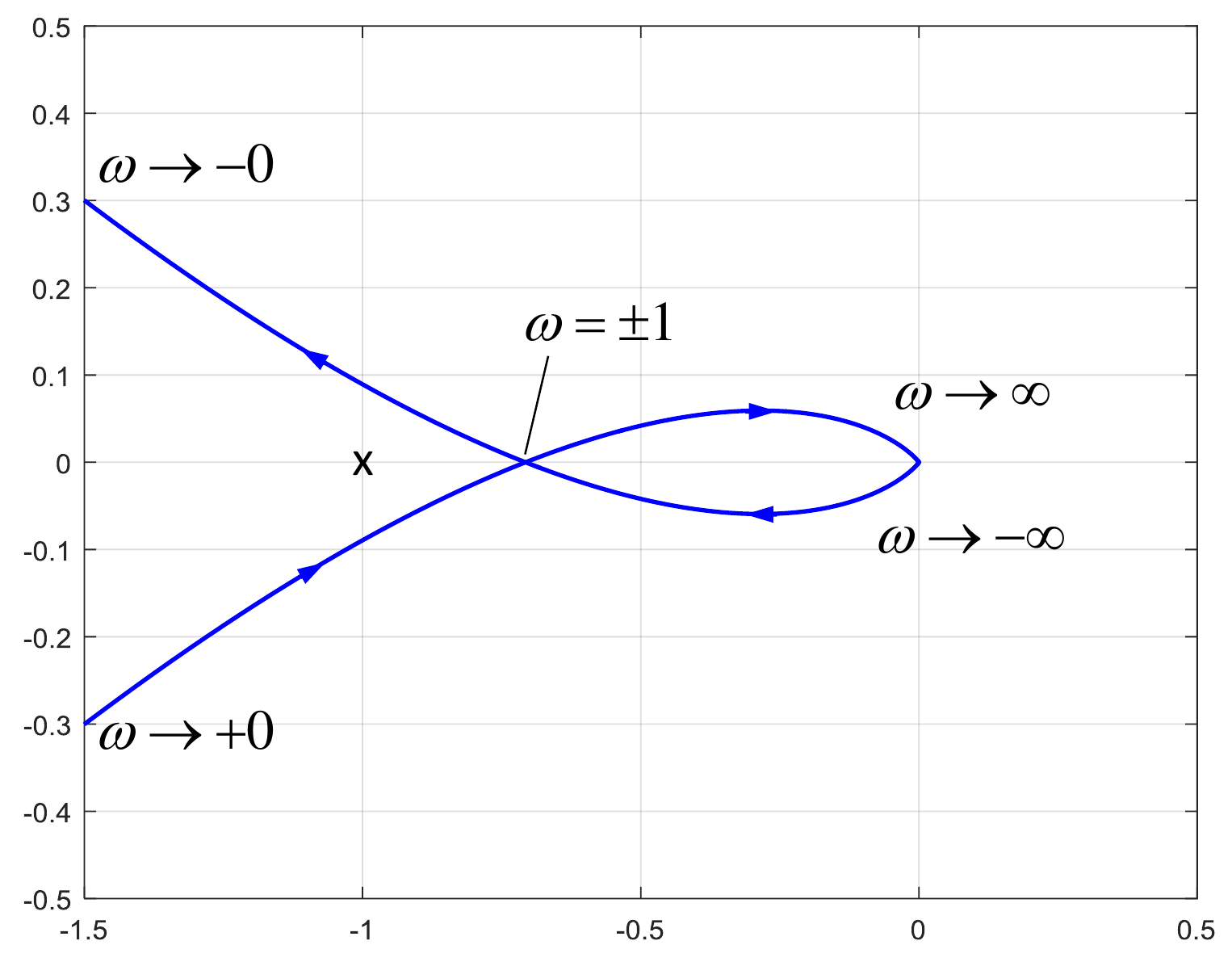}
	\caption{Nyquist plot of $z = F_o(j\omega) = \displaystyle \frac{K}{(j\omega)^{3/2}\,(j\omega+1)}$ (for $K=1$)}
	\label{fig:fig_example}
\end{figure}

For illustration, we consider an example where the open-loop transfer function is
\[ F_o(s) = \frac{K}{s^{3/2}\,(s+1)} \quad \mbox{with } K>0 \,. \] 
$F_o(s)$ has a pole in $s=-1$ and a polar branch point in $s=0$ with $\kappa_0 = -3/2$, thus $P_{o,+} = P_{o,j} =0$, $B_{o,j} = 3/2$. The Nyquist plot is shown in fig. \ref{fig:fig_example}; it intersects the negative real axis in the point $z = -K/\sqrt{2}$ for $\omega = \pm 1$. 
For $\omega \rightarrow 0$, $-1+F_o(j\omega) \approx K (j\omega)^{-3/2}$, thus $\phase{-1+F_o(j\omega)} \approx \phase{j^{-3/2}} = \pm 3/4 \pi$ ($-3/4 \pi$ for $\omega \rightarrow +0$ and $3/4 \pi$ for $\omega \rightarrow -0$).

We have to differentiate between 3 cases: $0<K<\sqrt{2}$, where the intersection point lies between $-1$ and $0$, $K=\sqrt{2}$ and $K>\sqrt{2}$. For $0<K<\sqrt{2}$, the change of argument is $\stackrel{+\infty}{\underset{\omega=-\infty}{\Delta^\bullet}} \phase{ 1 + F_o (j\omega)} = 3/2 \,\pi$ which coincides with $\pi \, \bigl( 2 P_{o,+}  + P_{o,j}  + B_{o,j} \bigr) = \pi \, B_{o,j}$, thus the closed-loop system is stable. In the other two cases, the change of argument is $-\pi/2$ (for $K=\sqrt{2}$) and $-5/2 \,\pi$ (for $K>\sqrt{2}$) and thus the closed-loop system is unstable.

\section{Conclusion}
\label{sec:Conclusion}
In this contribution, we considered linear, time-invariant SISO systems with non-rational transfer function from a perspective of complex analysis and singularities of the transfer function. We derived generalized criteria for BIBO stability and for the Nyquist test for closed-loop stability in terms of singularities with special focus on branch points on the imaginary axis. The class of transfer functions under consideration comprises systems with dead-time, retarded quasi-polynomial (RQ) meromorphic functions, fractional-order systems and the IO behavior of many systems with spatially distributed parameters. In context of BIBO stability, our approach makes a distinction between commensurate and incommensurate orders of branch points obsolete. Some limiting assumptions were made: for example we did not yet consider systems with impulse response containing distributions like in \cite{curtain_95}, and we excluded certain kinds of branch points like logarithms. This will be the focus of future work as well as a comparison of our approach to those of \cite{petras2008stability} and \cite{matignon1998stability}.

\section{Appendix}
\label{sec:Appendix}

\textbf{Proof of relation \eqref{eq:exp-stab}}

"$\Rightarrow$": By assumption, the system is $\beta$-exponentially stable, thus its impulse response satisfies $\abs{g(t)} \le M e^{-a t}$ where $M > 0$ and $a > - \beta \ge 0$ are known numbers. Choose $\beta_1$ between $-a$ and $\beta$: $-a < \beta_1 < \beta$. Thus
\[
\int \limits_{t=0}^\infty \abs{g(t)} e^{-\beta_1 t} dt 
\le
\int \limits_{t=0}^\infty M e^{-a t} e^{-\beta_1 t} dt 
=
M \int \limits_{t=0}^\infty e^{-(a+\beta_1) t} dt < \infty 
\]
since $a+\beta_1 > 0$. 

"$\Leftarrow$": Now by assumption, for given $\beta \le 0$, there exists $\beta_1 < \beta$ satisfying
\[
\int \limits_{t=0}^\infty \abs{g(t)} e^{-\beta_1 t} dt < \infty \, .
\]
This implies $\abs{g(t)} e^{-\beta_1 t} \rightarrow 0$ as $t\rightarrow\infty$. Since $g(t)$ is bounded by assumption, $\abs{g(t)} \le M e^{\beta_1 t} = M e^{-\alpha t}$ if we choose $a = - \beta_1$.

\xqed{$\blacksquare$}

\noindent \textbf{Proof of the increase of argument of imaginary singular points}

The Nyquist criterion derived in section \ref{sec:NyquistCrit}, makes use of the increase of argument when integrating along the semi-circles as shown in fig. \ref{fig:fig_Nyquist-contour} c). When calculating these values we have to differentiate between poles and the 4 kinds of branch points as already considered at the end of section \ref{sec:Closedloop}.

On the semi-circle  with center $a$ and radius $\varrho$ shown in fig. \ref{fig:fig_Nyquist-contour} c), $s = a + \varrho \, e^{j\,\varphi}$ and $ds = j \, \varrho \, e^{j\,\varphi} d\varphi = j\,(s-a)d\varphi$ with $-\pi/2 \le \varphi \le \pi/2$.

If $s=a$ is a pole of $F_o(s)$ with multiplicity $k>0$, then near $a$, 
\[F_o(s)=\frac{r_k}{(s-a)^k} \, \bigl(1+o(1)\bigr) \quad \mbox{and} \quad
F_o'(s) =-\frac{k\,r_k}{(s-a)^{k+1}} \, \bigl(1+o(1)\bigr) \, .\]
Thus, the integral along the semi-circle yields
\[
\int \limits_{C_\varrho} \frac{F_o'(s)}{1+F_o(s)}ds = 
\int \limits_{\varphi=-\pi/2}^{\pi/2} 
\frac{ \displaystyle
-\frac{k\,r_k}{(s-a)^{k+1}}\, \bigl(1+o(1)\bigr)}
{ \displaystyle \frac{r_k}{(s-a)^k}\, \bigl(1+o(1)\bigr)} j\, (s-a)  \, d\varphi
\rightarrow -j\,k\,\pi \,.
\]
Now, if $s=a$ is a zero of $1+F_o(s)$ (i.e. pole of $F_{cl}$) with multiplicity $k>0$, one obtains in a similar manner
\[
\int \limits_{C_\varrho} \frac{F_o'(s)}{1+F_o(s)}ds \rightarrow + j\,k\,\pi \,.
\]
If $s=a$ is a polar branch point of $F_o(s)$ with dominant exponent $\kappa_0<0$, then near $a$, 
\[F_o(s)=c_0\,(s-a)^{\kappa_0} \, \bigl(1+o(1)\bigr) \quad \mbox{and} \quad
F_o'(s) =\kappa_0 \,c_0 \, (s-a)^{\kappa_0-1} \, \bigl(1+o(1)\bigr) \, .\]
Thus, the integral along the semi-circle yields
\[
\int \limits_{C_\varrho} \frac{F_o'(s)}{1+F_o(s)}ds = 
\int \limits_{\varphi=-\pi/2}^{\pi/2} 
\frac{ \displaystyle
\kappa_0 \,c_0 \, (s-a)^{\kappa_0-1} \, \bigl(1+o(1)\bigr)}
{ \displaystyle c_0\,(s-a)^{\kappa_0}\, \bigl(1+o(1)\bigr)} j\, (s-a)  \, d\varphi
\rightarrow j\,\kappa_0\,\pi \,.
\]
If $F_o(s)$ has a regular branch point in $s=a$ with $\kappa_0 = 0$ and $c_0 \ne -1$,then 
\[F_o(s)=c_0+c_1\,(s-a)^{\kappa_1} \, \bigl(1+o(1)\bigr) \quad \mbox{and} \quad
F_o'(s) =\kappa_1 \,c_1 \, (s-a)^{\kappa_1-1} \, \bigl(1+o(1)\bigr) \, \]
and the integrand reads 
\begin{align}
\frac{F_o'(s)}{1+F_o(s)} j(s-a) &= 
\frac{ \displaystyle \kappa_1 \,c_1 \, (s-a)^{\kappa_1-1} \, \bigl(1+o(1)\bigr)}
{ \displaystyle 1+c_0+c_1\,(s-a)^{\kappa_1} \, \bigl(1+o(1)\bigr)} j\, (s-a)  \nonumber \\
&= \frac{\kappa_1 \, c_1}{1+c_0} (s-a)^{\kappa_1} \, j \,\bigl(1+o(1)\bigr) \, . \nonumber
\end{align}
Since $\kappa_1 > \kappa_0 = 0$, it tends to zero for $s\rightarrow a$. For $\kappa_0 = 0$ and $c_0 = -1$, we obtain
\[
\int \limits_{C_\varrho} \frac{F_o'(s)}{1+F_o(s)}ds = 
\int \limits_{\varphi=-\pi/2}^{\pi/2} 
\frac{ \displaystyle
\kappa_1 \,c_1 \, (s-a)^{\kappa_1-1} \, \bigl(1+o(1)\bigr)}
{ \displaystyle c_1\,(s-a)^{\kappa_1}\, \bigl(1+o(1)\bigr)} j\, (s-a)  \, d\varphi
\rightarrow j\,\kappa_1\,\pi \, ,
\]
and finally for $\kappa_0 > 0$:
\[F_o(s)=c_0\,(s-a)^{\kappa_0} \, \bigl(1+o(1)\bigr) \quad \mbox{and} \quad
F_o'(s) =\kappa_0 \,c_0 \, (s-a)^{\kappa_0-1} \, \bigl(1+o(1)\bigr) \, \]
thus the integral yields
\begin{align}
\int \limits_{C_\varrho} \frac{F_o'(s)}{1+F_o(s)}ds &= 
\int \limits_{\varphi=-\pi/2}^{\pi/2} 
\frac{ \displaystyle
\kappa_0 \,c_0 \, (s-a)^{\kappa_0-1} \, \bigl(1+o(1)\bigr)}
{ \displaystyle 1 + c_0\,(s-a)^{\kappa_0}\, \bigl(1+o(1)\bigr)} j\, (s-a)  \, d\varphi
\nonumber \\
&= \int \limits_{\varphi=-\pi/2}^{\pi/2} 
\kappa_0 \, c_0 \, \varrho^{\kappa_0}\, e^{j\,\kappa_0\,\varphi} \, j \,\bigl(1+o(1)\bigr) \, j \, d\varphi \; \rightarrow \, 0 \quad \mbox{as } \varrho \rightarrow 0 \, . \nonumber
\end{align}

\xqed{$\blacksquare$}

\bibliographystyle{alpha}
\bibliography{Bib}


\end{document}